\documentclass[12pt]{article} 
\usepackage{graphicx,psfrag,amsmath}

\def\N{I\!\!N}

\newcommand{\qed}{\nobreak \ifvmode \relax \else
\ifdim\lastskip<1.5em \hskip-\lastskip
\hskip1.5em plus0em minus0.5em \fi \nobreak
\vrule height0.75em width0.5em depth0.25em\fi}

\begin{document} 

\vspace*{-2.8cm}

\begin{center}
{\Large {\bf {On The Discrepancy of Quasi-progressions}}}

\vskip 15 pt

{\large {\em {Sujith Vijay}}}

\vskip 15 pt

{\large {\tt {Department of Mathematics}}},

\vskip 10 pt

{\large {\tt {Rutgers, the State University of New Jersey}}}
\end{center}

\vskip 10 pt

A classic result in discrepancy theory is the theorem of Roth [2] stating
that if the elements of $\{0,1,2,\ldots, n \}$ are $2$-coloured, there
exists an arithmetic progression $\{a, a+d, \ldots, a+(k-1)d\}$ of
discrepancy at least $(1/20) n^{1/4}$, with $0 \le a < d \le \sqrt{6n}$.
In 1996, Matou\v{s}ek and Spencer [8], building upon results by
S\'ark\"ozy (see [3]) and Beck [5], showed that apart from constants, this
result is the best possible. \\

The situation is quite different, however, for homogeneous arithmetic
progressions (HAPs), the subfamily of arithmetic progressions containing
$0$. It turns out that there are extremely balanced colourings for such
arithmetic progressions. Consider the following examples.

$$\chi_3(3k+1)=1; \; \chi_3(3k+2)=-1; \; \chi_3(3k)= \chi_3(k)$$
$$\chi^{*}_3(3k+1)=1; \; \chi^{*}_3(3k+2)=-1; \; \chi^{*}_3(3k)= 0 \quad 
\; \;$$

It is easy to show that all HAPs contained in $\{0,1,\ldots, n\}$ have
discrepancy $O(\log \, n)$ under $\chi_3$ and discrepancy at most $1$
under $\chi^{*}_3$. (Of course, $\chi^{*}_3$ is a $3$-colouring, and if we
are going to allow $0$ we might as well colour everything $0$, but let us
look the other way and at the bright side for a moment.) In general, for
any prime $p$, we can define a colouring $\chi^{*}_p$ via the non-trivial
real character modulo $p$ (the Legendre symbol). The discrepancy of all
HAPs are bounded by $(p-1)/2$ under this colouring.  Also note that
$\chi^{*}_p$ is ``almost admissible" for large $p$, since only a small
fraction of numbers is coloured $0$. Whether there is a ``completely
admissible" colouring of bounded discrepancy for HAPs is a question raised
by Erd\H{o}s in the 1930s, and one that remains unsolved to this day. It
is indeed a mishap that this innocent-looking question should turn out to
be so difficult. \\

Yet homogeneous arithmetic progressions are tiny herrings on the tip of
the iceberg of quasi-progressions. Perhaps a definition is in order. A
quasi-progression $Q(\alpha;s,t)$ is the sequence of integers $$\lfloor s
\alpha \rfloor, \lfloor (s+1) \alpha \rfloor, \ldots, \lfloor t \alpha
\rfloor$$

In other words, a quasi-progression is a sequence of successive multiples
of a real number, with each multiple rounded down to the nearest integer.
Since distinct real numbers yield distinct quasi-progressions, we are
dealing with an uncountable family of sequences. Note that for integer
values of $\alpha$, quasi-progressions reduce to HAPs, or the
set-difference of two HAPs. Thus the problem raised by Erd\H{o}s concerns
a subfamily of quasi-progressions, corresponding to integer values of
$\alpha$. \\

Our first theorem gives a lower bound on the discrepancy of the family of
all quasi-progressions contained in $\{0,1,\ldots,n\}$. \\

{\bf {Theorem 1}} If the integers from $0$ to $n$ are $2$-coloured, there
exists $\alpha > 1$ and integers $s$ and $t$ such that the
quasi-progression $Q(\alpha;s,t)$ has discrepancy at least $(1/50)
n^{1/6}$. \\

{\bf {Proof}} Let $m < n$. The value of $m$ will be specified at the end
of the proof. By Roth's theorem, there exists an arithmetic progression
$P_1=\{ a,a+d, a+2d, \ldots \}$ contained in $\{0,1,\ldots,m \}, \, 2 \le
d < \sqrt{6m}, \, 0 \le a < d$, with discrepancy at least $(1/40)  
m^{1/4}$.  Let $P_2=(n-m)+P_1$. We will show that for suitably chosen $m$,
$P_2$ can be realised as a quasi-progression corresponding to a real
number $\alpha > 1$. \\

Observe that if $\alpha = d - \epsilon$, the first $\lfloor 1/ \epsilon
\rfloor$ elements in the sequence $\lfloor \alpha \rfloor, \lfloor 2
\alpha \rfloor, \lfloor 3 \alpha \rfloor \ldots$ are congruent to $-1
(\mbox{mod }d)$, the next $(\lfloor 2/ \epsilon \rfloor - \lfloor 1/
\epsilon \rfloor)$ elements are congruent to $-2 \, (\mbox{mod }d)$, and
so on. In particular, the arithmetic progression $P_2 \equiv - (d-a) \,
(\mbox{mod }d)$ can be realised as a quasi-progression by choosing
$\epsilon$ such that $P_2$ is completely contained in the $(d-a)^{th}$
block of length $(1/\epsilon) + O(1)$. \\

Since $P_2 \subseteq \{n-m, n-m+1, \ldots, n \}$, it suffices to choose
$\epsilon$ such that $(d-a-1) \lceil 1/\epsilon \rceil < n-m$ and $(d-a)
\lfloor 1/ \epsilon \rfloor > n$. Such an $\epsilon$ exists if
$$\frac{n-m}{d-a-1}-\frac{n}{d-a} > 1$$

Note that $d-a \le d \le \sqrt{6m}$. Therefore, we can choose $m=\lfloor
6^{-1/3}n^{2/3} \rfloor$. This yields a quasi-progression of discrepancy at 
least $(1/50)n^{1/6}$. \qed \\

While it is not known whether the set of homogeneous arithmetic
progressions have bounded discrepancy, there exist colourings (see [9])
for which the arithmetic progression $\{0,d,2d,\ldots,\}$ has
discrepancy at most $d^{4+o(1)}$ for all $d$. It turns out, however, that
upper bounds independent of $n$ do not exist for most quasi-progressions.
\\

Let $\alpha > 1$ be given, together with a $2$-colouring of
$\{0,1,\ldots,n\}$. Let $D_{\alpha}(n)$ denote the maximum discrepancy of
$Q(\alpha;s,t)$ over all admissible $s$ and $t$. 
In 1986, Beck [6] showed that given
any $2$-colouring of the non-negative integers, for almost every $\alpha
\in [1, \infty)$, there are infinitely many $n$ such that $D_{\alpha}(n)
\ge \log^{*} n$. Recall that $\log^{*} x $ denotes the inverse of the
tower function:  $\log^{*} x = \ln x$ for $1 < x < e$ and $\log^{*} (e^x)
= 1 + \log^{*} x$. \\

We improve on this result, and prove the following theorem. \\

{\bf {Theorem 2}} Let $\chi$ be a partial colouring of the non-negative
integers with density $\rho > 0$, and let $\chi_n$ denote the restriction
of $\chi$ to $\{0,1,\ldots, n\}$. Then for almost every $\alpha \in [1,
\infty)$, there are infinitely many $n$ such that $D_{\alpha}(n) \ge (\log
n)^{1/4-o(1)}$. \\

{\bf {Proof}} Let $E$ denote the set of $\alpha$ such that there are only
finitely many $n$ with $D_n(\alpha) \le (\log n)^{1/4-o(1)}$ under the
colouring $\chi$. If $E$ has positive measure, there exists a positive
integer $t$ for which the set of balanced $\alpha$ in $[t,t+1)$ has
measure $\delta > 0$. But it follows from the Main Lemma (see below) that
there exists $c_0=c_0(\rho, t, \delta)$ such that the set of $\alpha$ with
$D_n(\alpha) \le c_0 (\log n)^{1/4}$ has measure less than $\delta$. For
all other $\alpha$ in $[t,t+1)$, we have $D_n(\alpha) > (\log
n)^{1/4-o(1)}$ for sufficiently large $n$, yielding a contradiction. \qed
\\

It remains to state and prove the Main Lemma. \\

{\bf {Main Lemma}} Let $\chi$ and $\chi_n$ be as in the statement of
Theorem 2. Given $t \in [1, \infty)$ and $\delta > 0$, there exists
$c_0=c_0(\rho, t, \delta)$ such that the set of $\alpha$ in $[t,t+1)$ with
$D_{\alpha}(n) \le c_0(\log n)^{1/4}$ under $\chi_n$ has Lebesgue measure
less than $\delta$. \\

{\em {Remark:}} We say that $\alpha$ is $M$-balanced if $D_{\alpha}(n) \le
M$.  For brevity, we shall hereafter refer to $(c_0(\log
n)^{1/4})$-balanced $\alpha$ simply as ``balanced". We will transform the
problem into a geometric setting, with a view to using orthogonal
functions, as was done by Roth [1] in his classic paper on the
measure-theoretic discrepancy of axis-parallel rectangles. A similar
construction was used by Hochberg [7] to show the existence of a
quasi-progression of discrepancy $c'_0 (\log n)^{1/4}$. As we saw in
Theorem 1, quasi-progressions with much larger discrepancy do occur. \\

{\bf {Proof}} We shall assume, for the sake of convenience, that
$n=(t+1)m$ where $m=2^u$ for some positive integer $u$. We join each
lattice point $(a,b)$ with the one vertically above it, and give the
resulting unit segment the colour $\chi(b)$. For each point $(x,y)$ in the
plane, the discrepancy function $D(x,y)$ is defined to be the sum of the
$\chi$-values of the unit segments crossed by the line joining $(0,0)$ and
$(x,y)$. Note that $|D(x,y)| \le M$ if and only if $y/x$ is $M$-balanced.
\\

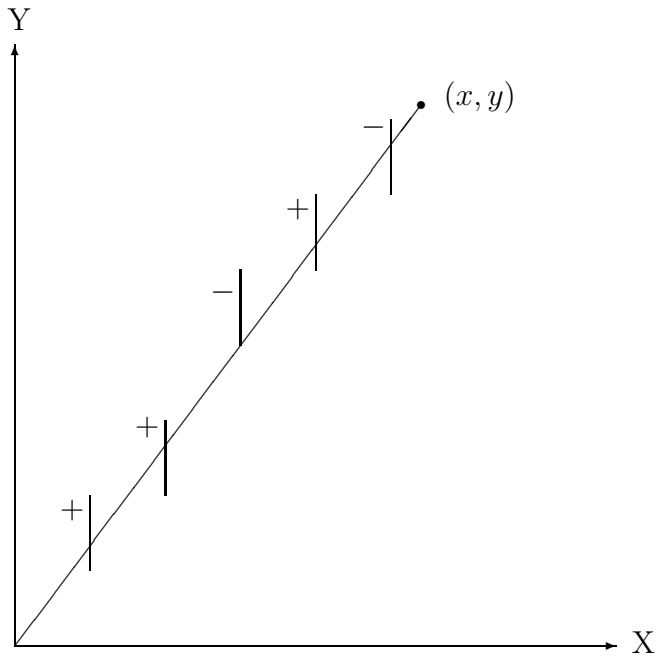
\setlength{\unitlength}{1mm}
\begin{picture}(90,100)(-10,-10)
\put(0,0){\vector(1,0){80}}
\put(0,0){\vector(0,1){80}}
\put(82,-1){\makebox{X}}
\put(-1,82){\makebox{Y}}
\put(0,0){\line(3,4){54}}
\put(54,72){\circle*{1}}
\put(57,72){\makebox{$(x,y)$}}
\put(10,10){\line(0,1){10}}
\put(20,20){\line(0,1){10}}
\put(30,40){\line(0,1){10}}
\put(40,50){\line(0,1){10}}
\put(50,60){\line(0,1){10}}
\put(6,17){\makebox{$+$}}
\put(16,28){\makebox{$+$}}
\put(26,46){\makebox{$-$}}
\put(36,57){\makebox{$+$}}
\put(46,68){\makebox{$-$}}
\end{picture}

Let $H(x,y) = D(x,y)$ if $y/x$ is balanced, and $0$ otherwise. Suppose
that the measure of the set of balanced $\alpha$ in $[t,t+1)$ is at least
$\delta$. We will deduce a contradiction for a suitably chosen $c_0$ by
producing a point $(x_0,y_0)$ with $H(x_0,y_0) > c_0 (\log n)^{1/4}$. \\

Let $R$ denote the region bounded by the lines $x=m/2, x=m, y=tx,
y=(t+1)x$. We will construct orthonormal functions $g_1, g_2, \ldots, g_r$
on $R$ where $r=(\log n)/8$ and $$\sum_{i=1}^r (\langle H,g_i \rangle)^2
\ge \frac{\rho^2 \delta^{13} m^2 (\log n)^{1/2}}{2^{27} c^2_0 t^3}$$

Since $R$ has area $3m^2/8$, it follows from Bessel's inequality
that there exists $(x_0,y_0)$ with $H(x_0,y_0) > c_0 (\log n)^{1/4}$ for
$$c_0(\rho,t,\delta)= \frac{{\rho^{1/2}} \delta^{13/4}}{142 t^{3/4}}$$
yielding the desired contradiction. \\

The functions $g_1, g_2, \ldots, g_r$ will be normalised versions of
mutually orthogonal functions $G_1, G_2, \ldots, G_r$. Following Hochberg,
we will construct $G_i$ by dividing $R$ into a grid of trapezoids, called
the $i^{th}$ trapezoidal grid. We use vertical lines spaced $\ell = 2^i$
apart and slanting lines with slopes equally spaced between $t$ and $t+1$.
The slopes of consecutive slanting lines differ by $\tau \doteq 1 / (\ell
\beta m)$ where $\beta = c_2 (\log n)^{1/4}$. The value of $c_2$ will be
specified later. It is easy to see that the individual grid trapezoids
have area at most $1 / \beta$ and at least $1/(2 \beta)$. \\

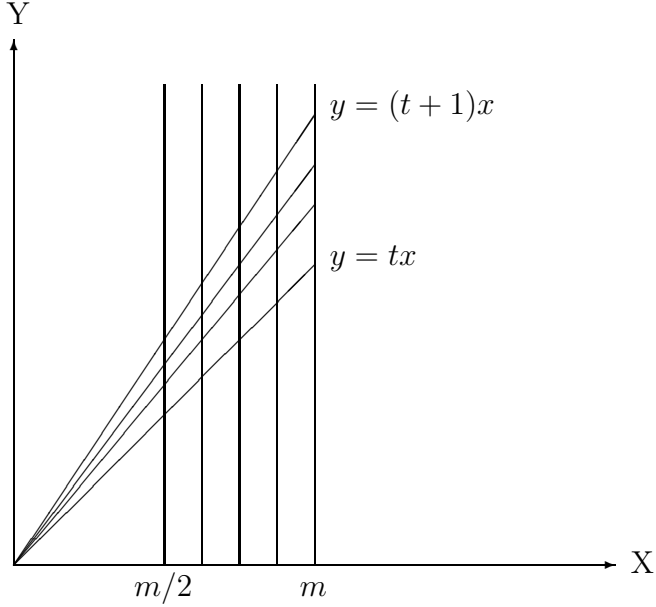
\setlength{\unitlength}{1mm}
\begin{picture}(90,80)(-10,-8)
\put(0,0){\vector(1,0){80}}
\put(0,0){\vector(0,1){70}}
\put(82,-1){\makebox{X}}
\put(-1,72){\makebox{Y}}
\put(0,0){\line(2,3){40}}
\put(0,0){\line(3,4){40}}
\put(0,0){\line(5,6){40}}
\put(0,0){\line(1,1){40}}
\put(42,40){\makebox{$y=tx$}}
\put(42,60){\makebox{$y=(t+1)x$}}
\put(16,-4){\makebox{$m/2$}}
\put(38,-4){\makebox{$m$}}
\multiput(20,0)(5,0){5}{\line(0,1){64}}
\end{picture}

Note that we have specified only the spacing between the grid lines and
not their actual position. We choose the position of the rightmost
vertical line randomly and uniformly in the interval $[n-\ell,n)$, and the
slope of the lowermost line randomly and uniformly in the interval $[t,t+
\tau)$. The region between two consecutive sloping lines will be called a
{\em {sector}}, and sectors will be identified with subintervals of
$[t,t+1)$ in the natural fashion. We will denote the measure of balanced
$\alpha$ in the $j^{th}$ sector of the $i^{th}$ grid by $\mu_{ij}$. For
convenience, we define $\mu^{*}_{ij} = \mu_{ij} / \tau$. \\

If $\chi(b) \neq \chi(b-1)$, we refer to $b$ as a {\em {switch value}}.  
Furthermore, a lattice point $(a,b)$ will be called a {\em {switch point}}
if $b$ is a switch value. A switch point is said to be {\em {good}} if it
finds itself alone in a trapezoid no matter how the grid is positioned;  
{\em {bad}} otherwise. We shall denote the number of good switch points in
the $j^{th}$ sector of the $i^{th}$ grid by $s^{*}_{ij}$. \\

We define $G_i$ as follows: On a trapezoid containing exactly one switch 
point, $G_i$ is defined in a checkerboard fashion. On all other 
trapezoids, $G_i$ is defined to be identically zero. \\

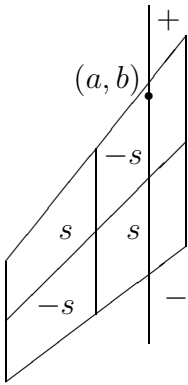
\setlength{\unitlength}{1mm}
\begin{picture}(86,52)(-40,15)
\put(10,20){\line(0,1){16}}
\put(22,29){\line(0,1){22}}
\put(34,38){\line(0,1){28}}
\put(10,20){\line(4,3){24}}
\put(10,28){\line(1,1){24}}
\put(10,36){\line(4,5){24}}
\put(26,39){\makebox{$s$}}
\put(23,49){\makebox{$-s$}}
\put(17,39){\makebox{$s$}}
\put(14,29){\makebox{$-s$}}
\put(29,25){\line(0,1){45}}
\put(31,30){\makebox{$-$}}
\put(30,67){\makebox{$+$}}
\put(29,58){\circle*{1}}
\put(19,59){\makebox{$(a,b)$}}

\end{picture}

The vertical dividing line passes through the centre of the trapezoid. The
position of the slanting dividing line is chosen such that the measure of
balanced $\alpha$ above the line and inside the sector equals the measure
of balanced $\alpha$ below the line and inside the sector. The value of
$s$ will vary from trapezoid to trapezoid, but will always equal $+1$ or
$-1$. Since the vertical dividing lines are nested dyadically (note that
the vertical spacing is $\ell=2^i$), it is clear that $\{G_i\}$ form an
orthogonal family. \\

We now derive a lower bound on the inner product $\langle H,G_i \rangle$.  
The position of the slanting dividing line has been chosen with a view to
extending Hochberg's argument for the $\mu^{*}_{ij}=1$ case to the more
general problem at hand. \\

\noindent{\bf{Lemma 1}} $E(\langle H,G_i \rangle) \ge (\underset{j}{\sum}
(\mu^{*}_{ij})^2 s^{*}_{ij}) /(32 \beta)$ \\

\noindent{\bf{Proof}} Consider the contribution of a unit vertical segment
$\ell_{a,b}$ joining $(a,b)$ and $(a,b+1)$ to the discrepancy function
$H(x,y)$. Let $$B_{a,b} = \left \{ (x,y):x \ge a, \frac{b}{a} \le
\frac{y}{x} < \frac{b+1}{a} \right \}$$ denote the set of points behind
the line $\ell_{a,b}$. \\

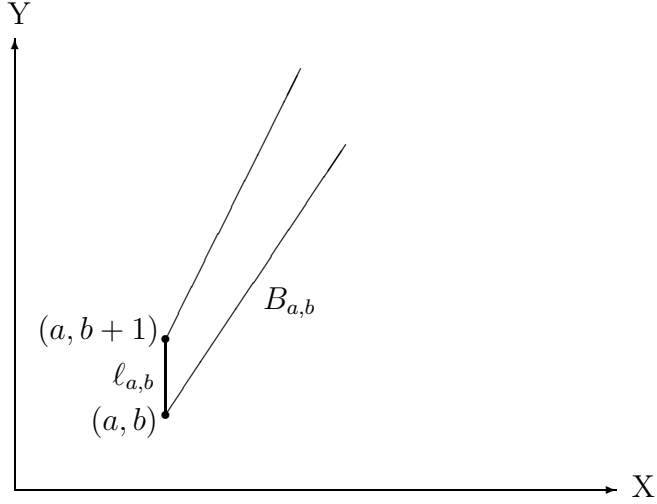
\setlength{\unitlength}{1mm}
\begin{picture}(100,75)(-10,-10)
\put(0,0){\vector(1,0){80}}
\put(0,0){\vector(0,1){60}}
\put(82,-1){\makebox{X}}
\put(-1,62){\makebox{Y}}
\put(20,10){\line(0,1){10}}
\put(20,10){\line(2,3){24}}
\put(20,20){\line(1,2){18}}
\put(20,10){\circle*{1}}
\put(20,20){\circle*{1}}
\put(10,8){\makebox{$(a,b)$}}
\put(3,20){\makebox{$(a,b+1)$}}
\put(13,14){\makebox{$\ell_{a,b}$}}
\put(33,24){\makebox{$B_{a,b}$}}

\end{picture}

Now define $$H_{a,b}(x,y)= \left\{ \begin{array}{ll} \chi(b) & \quad
\mbox{if $(x,y) \in B_{a,b}$ and $y/x$ is balanced} \\ 0 & \quad
\mbox{otherwise} \\ \end{array} \right.$$ 

Clearly, $$H(x,y) = \sum_{a=0}^{\infty} \sum_{b=0}^{\infty} H_{a,b}(x,y)$$
Furthermore, only finitely many terms in this sum are non-zero, for any
fixed $(x,y)$. Consider a good switch point $(a,b)$ lying inside a
trapezoid $T$, located in the $j^{th}$ sector of the $i^{th}$ grid. \\

We claim that if neither $(a,b)$ nor $(a,b+1)$ lie inside $T$, then $$
\underset{T \; \; } {\int \int} G_i(x,y) H_{a,b}(x,y) \, dx \, dy = 0$$

If $T$ lies entirely outside or entirely inside $B_{a,b}$, it is clear
that the integral is zero. If exactly one of the bounding lines of
$B_{a,b}$ intersects $T$, the geometric symmetry with respect to the 
vertical dividing line or the measure-theoretic symmetry with respect to 
the sloping dividing line, as the case may be, ensures that there is 
perfect cancellation. Thus the integral vanishes in this case as well. \\

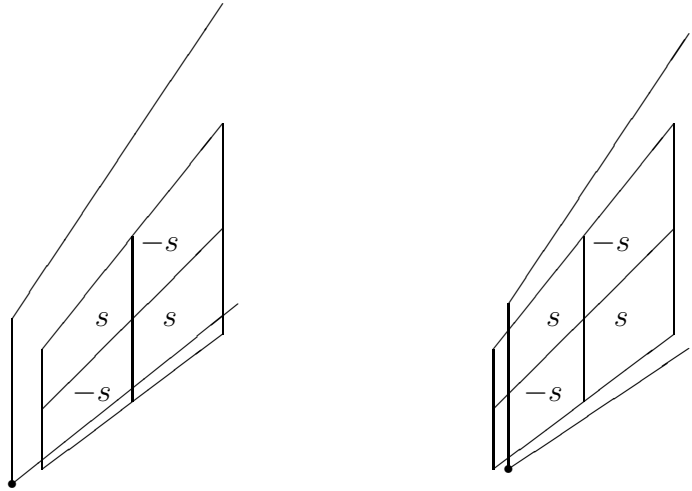
\begin{picture}(120,65)(-10,15)
\put(10,20){\line(0,1){16}}
\put(22,29){\line(0,1){22}}
\put(34,38){\line(0,1){28}}
\put(10,20){\line(4,3){24}}
\put(10,28){\line(1,1){24}}
\put(10,36){\line(4,5){24}}
\put(26,39){\makebox{$s$}}
\put(23,49){\makebox{$-s$}}
\put(17,39){\makebox{$s$}}
\put(14,29){\makebox{$-s$}}
\put(70,20){\line(0,1){16}}
\put(82,29){\line(0,1){22}}
\put(94,38){\line(0,1){28}}
\put(70,20){\line(4,3){24}}
\put(70,28){\line(1,1){24}}
\put(70,36){\line(4,5){24}}
\put(86,39){\makebox{$s$}}
\put(83,49){\makebox{$-s$}}
\put(77,39){\makebox{$s$}}
\put(74,29){\makebox{$-s$}}
\put(6,18){\circle*{1}}
\put(6,18){\line(0,1){22}}
\put(6,18){\line(5,4){30}}
\put(6,40){\line(2,3){28}}

\put(72,20){\circle*{1}}
\put(72,20){\line(0,1){22}}
\put(72,20){\line(3,2){24}}
\put(72,42){\line(2,3){24}}
\end{picture}

Therefore, we need consider only the terms $H_{a,b}(x,y)$ and
$H_{a,b-1}(x,y)$, where $(a,b) \in T$. If $(a,b)$ is not a switch point,
we have, $$\underset{T \; \; } {\int \int} G_i(x,y) (H_{a,b}(x,y) +
H_{a,b-1}(x,y))  \, dx \, dy = 0$$

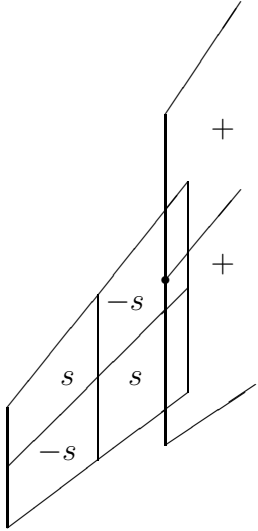
\begin{picture}(120,70)(-30,16)
\put(10,20){\line(0,1){16}}
\put(22,29){\line(0,1){22}}
\put(34,38){\line(0,1){28}}
\put(10,20){\line(4,3){24}}
\put(10,28){\line(1,1){24}}
\put(10,36){\line(4,5){24}}
\put(26,39){\makebox{$s$}}
\put(23,49){\makebox{$-s$}}
\put(17,39){\makebox{$s$}}
\put(14,29){\makebox{$-s$}}
\put(31,53){\circle*{1}}
\put(31,31){\line(0,1){44}}
\put(31,53){\line(5,6){10}}
\put(31,31){\line(3,2){12}}
\put(31,75){\line(2,3){10}}
\put(37,72){\makebox{$+$}}
\put(37,54){\makebox{$+$}}
\end{picture}

Now suppose that $(a,b)$ is a switch point. If $(a,b)$ lies on the
intersection of the two dividing lines, we have

$$\underset{T \; \; } {\int \int} G_i(x,y) (H_{a,b}(x,y) + H_{a,b-1}(x,y))  
\, dx \, dy = \frac{s}{4}(\chi(b)-\chi(b-1))\mu^{*}_{ij} \mbox{ area}(T)
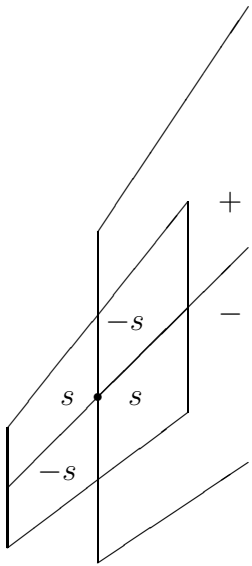
$$

\begin{picture}(86,72)(-30,16)
\put(10,20){\line(0,1){16}}
\put(22,29){\line(0,1){22}}
\put(34,38){\line(0,1){28}}
\put(10,20){\line(4,3){24}}
\put(10,28){\line(1,1){24}}
\put(10,36){\line(4,5){24}}
\put(26,39){\makebox{$s$}}
\put(23,49){\makebox{$-s$}}
\put(17,39){\makebox{$s$}}
\put(14,29){\makebox{$-s$}}
\put(22,40){\circle*{1}}
\put(22,18){\line(0,1){44}}
\put(22,40){\line(1,1){20}}
\put(22,18){\line(3,2){20}}
\put(22,62){\line(2,3){20}}
\put(38,65){\makebox{$+$}}
\put(38,50){\makebox{$-$}}
\end{picture}

We choose $s$ so that the integral is positive. Since the switch point $p$
is good, there are no other lattice points in $T$, and the value of $s$
can now be safely assumed fixed. Thus we get $$ \underset{T \; \; } {\int
\int} G_i(x,y) H(x,y)) \, dx \, dy \ge \frac{\mu^{*}_{ij}}{8 \beta} $$
provided $(a,b)$ lies on the intersection of the two dividing lines. Since
the location of $(a,b)$ inside the trapezoid is a uniformly distributed
random variable, we have $$E \left ( \underset{T \; \; } {\int \int}
G_i(x,y) H(x,y) \, dx \, dy \right ) \ge \frac{(\mu^{*}_{ij})^2}{32
\beta}$$

Adding over all switch points and using the linearity of expectation, we
get $$E( \langle H,G_i \rangle) \ge \frac{1}{32 \beta} \sum_j
(\mu^{*}_{ij})^2 s^{*}_{ij}$$ as claimed. \qed \\

We now prove a slightly stronger version of a lemma due to Beck [6]. \\

\noindent {\bf {Lemma 2}} Let $J \subseteq [0,1]$ be an arbitrary interval
of length $\lambda$ and let $1 \le b_1 < b_2 < \ldots b_q$ be integers.
Let $N(\alpha,J) = |\{j: \{b_j \alpha \} \in J, 1 \le j \le q \} |$. If $q
\ge \lambda^{-6}$, then $\mu(\{ \alpha \in [0,1]: N(\alpha,J)  \ge (q
\lambda /2) \}) \ge 1 - (8 / \sqrt{q})$. \\

\noindent {\bf {Proof}} The proof uses LeVeque's inequality from the
theory of uniform distributions, and is almost identical to the proof of
Beck's original lemma. \\

Let $x_j = \{b_j \alpha \}, 1 \le j \le q$. Define $\Delta(\alpha)$ and
$S_n(\alpha)$ as follows: $$\Delta(\alpha)=\sup_{0 \le a < b \le 1} \left
| \left ( \sum_{j: x_j \in [a,b)} \frac{1}{q} \right) - (b-a) \right |$$

$$S_n(\alpha) = \frac{1}{q} \sum_{j=1}^q e^{2 \pi i n x_j}$$

Note that $$\int_0^1 |S_n(\alpha)|^2 \, d \alpha = \frac{1}{q^2} \int_0^1
\sum_{j=1}^q \sum_{k=1}^q e^{2 \pi i n (b_j - b_k) \alpha} \, d \alpha =
\frac{1}{q}$$

By LeVeque's inequality, $$\Delta^3(\alpha) \le \frac{6}{\pi^2} \sum_{n
\in N} \frac{ |S_n(\alpha)|^2}{n^2}$$

Therefore, $$\int_0^1 \Delta^3(\alpha) \, d \alpha \le \frac{6}{\pi^2}
\int_0^1 \left ( \sum_{n \in \N} \frac{1}{n^2} |S_n(\alpha)|^2 \right ) \,
d \alpha = \frac{1}{q}$$

Let $E=\{\alpha \in [0,1): N(\alpha,J) \ge q \lambda / 2 \}$ and
$F=[0,1) \setminus E$. Clearly,

$$\frac{\lambda^3 \mu(F)}{8} \le \int_0^1 \left | \left ( \sum_{j:  x_j
\in J} \frac{1}{q} \right ) - \lambda \right | \, d \alpha \le \int_0^1
\Delta^3(\alpha) \, d \alpha$$

Therefore, $\lambda^3 \mu(F) / 8 \le 1/q$. Since $q \ge
\lambda^{-6}$, we have $$\mu(E) = 1 - \mu(F) \ge 1 - \frac{8}{\sqrt{q}},$$
proving the lemma. \qed \\

Let $b_1, b_2, \ldots, b_q$ be the switch values of the colouring $\chi$
in $[N/2,N]$. Note that $q \ge (N \rho)/ (4 c_0 (t+1)(\log n)^{1/4}) = (m
\rho)/(4 c_0 (\log n)^{1/4})$. Since switch points come in rows, it is
clear that $\chi_n$ gives rise to $mq$ switch points. \\

\noindent{\bf {Lemma 3}} $\underset{j}{\sum} (\mu^{*}_{ij})^2 s^{*}_{ij}
\ge \delta^5 mq/(4096 t)$, for $1 \le i \le r$ \\

\noindent{\bf {Proof}} We say that a sector is {\em {rich}} if
$\mu^{*}_{ij} > \delta / 2$.  Let $\{I_k\}_{k=1}^L$ be an enumeration of
the rich sectors. Since $(\delta/2) (n \ell \beta - L) + L \tau \ge \delta
\, ,$ we have $L > (\delta / 2) m \ell \beta$. \\

We use Lemma 2 with $J = [0, \delta/(4 \ell \beta (t+1)]$, so that
$\lambda = \delta / (4 \ell \beta (t+1))$. Since $r= \log m /8$ and $i \le
r$, we have $q \ge \lambda^{-6}$. For an arbitrary interval $I=[a_k,b_k)$,
let $I'$ and $I''$ denote $[a_k, c_k)$ and $[c_k,b_k)$ respectively, where
the measure of balanced $\alpha$ in $I'$ and $I''$ are equal. Let $A =
\cup_{k=1}^L I'_k$. Note that $A$ has measure at least
$\delta^2/8$. Let $B=\{ \theta: 1/ \theta \in A \}$. Since $A
\subseteq [t, t+1) \subseteq [t,2t)$, the measure of $B$ is at least
$\delta^2/(32t^2)$. \\

Let $B^{*} = \{ \theta \in B: N(\theta,J) \ge (q \lambda) / 2 \}$. For
sufficiently large $n$, $B^{*}$ has measure at least
$\delta^2/(64t^2)$. Note that $\theta \in B^{*} \Rightarrow 1/
\theta \in I'_k$ for some $k$. Suppose $\{b_v \theta\} \in J$ for such a
$\theta$. Let $a_v = \lfloor b_v \theta \rfloor$. Then we have, 
$$0 < \frac{b_v}{a_v} - \frac{1}{\theta} < \frac{2 \lambda (t+1)}{m} < 
\frac{\delta \tau}{4}$$

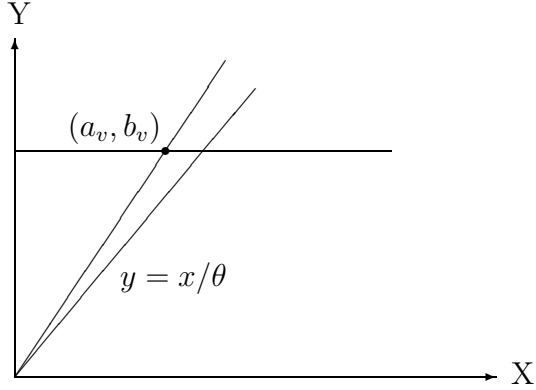
\begin{picture}(90,54)(-10,-5)
\put(0,0){\vector(1,0){64}}
\put(0,0){\vector(0,1){45}}
\put(66,-1){\makebox{X}}
\put(-1,47){\makebox{Y}}
\put(0,0){\line(2,3){28}}
\put(0,0){\line(5,6){32}}
\put(0,30){\line(1,0){50}}
\put(20,30){\circle*{1}}
\put(7,32){\makebox{$(a_v,b_v)$}}
\put(14,12){\makebox{$y=x / \theta$}}
\end{picture}

It follows that $(b_v/a_v) \in I_k$. Thus the $k^{th}$ sector contains a
switch point of the form $(a_v, b_v)$. \\

Since the contribution of a single $I'_k$ towards the measure of $B^{*}$
is at most $\tau/t^2$, there must be at least $\delta^2 m \ell
\beta/(32t^2)$ rich sectors contributing at least $\delta^3 mq/(1024t)$
switch points between them. \\

We now derive an upper bound on the total number of bad switch points. \\

Given a bad switch point $(a,b)$, there exists $a'$ such that $$|a' - a| 
\le \ell \mbox{ and } \left | \left | \frac{ba'}{a} \right | \right | < 
\frac{1}{\ell \beta}
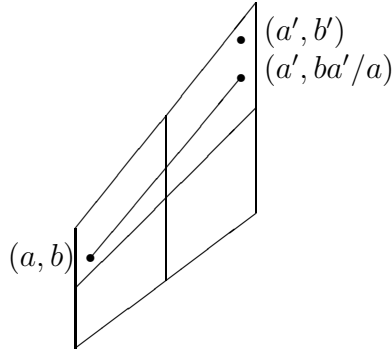
$$

\begin{picture}(86,52)(-30,20)
\put(10,20){\line(0,1){16}}
\put(22,29){\line(0,1){22}}
\put(34,38){\line(0,1){28}}
\put(10,20){\line(4,3){24}}
\put(10,28){\line(1,1){24}}
\put(10,36){\line(4,5){24}}
\put(12,32){\circle*{1}}
\put(32,56){\circle*{1}}
\put(32,61){\circle*{1}}
\put(35,56){\makebox{$(a',ba'/a)$}}
\put(35,61){\makebox{$(a',b')$}}
\put(12,32){\line(5,6){20}}
\put(1,31){\makebox{$(a,b)$}}
\end{picture}

Let $d=|a'-a|$. Note that there are $m \ell/2$ pairs $(a,d)$ with $1 \le
d \le \ell, \, m/2 \le a \le m$.  For each such pair, there are
$a/(\ell \beta) + O(\ell)$ values of $b$ that satisfy
$||bd/a|| < 1/(\ell \beta)$. It follows that there are at
most $m^2 / \beta$ lattice points which do not find themselves alone
in a trapezoid for some placement of the grid. \\

Since the number of bad switch points is at most $$\frac{m^2}{\beta} =
\frac{m^2}{c_2 (\log n)^{1/4}} < \frac{\delta^3 mq}{1024t}$$ for $c_2 \ge
4096c_0 t / (\delta^3 \rho)$, we have $$\sum (\mu^{*}_{ij})^2 s^{*}_{ij} >
\frac{\delta^5 mq}{4096t}$$ as required. \qed \\

Note that $$E({\langle H,G_i \rangle}^2) \ge [E(\langle H,G_i \rangle)]^2
\ge \left ( \sum_j \frac{(\mu^{*}_{ij})^2 s^{*}_{ij}}{32 \beta} \right
)^2$$

Furthermore, the combined area of all the grid trapezoids containing
exactly one switch point is at most $mq / \beta$. Therefore, $$E({\langle
H,g_i \rangle}^2) \ge \sum_j \frac{(\mu^{*}_{ij})^2 s^{*}_{ij})^2}{1024
mq \beta} \ge \frac{\rho^2 \delta^{13} m^2}{2^{24} c_0^2 t^3 (\log
n)^{1/2}}$$

By the linearity of expectation, $$E \left ( \sum_{i=1}^r {\langle H,g_i
\rangle}^2 \right ) \ge \frac{\rho^2 \delta^{13} m^2 (\log
n)^{1/2}}{2^{27} c_0^2 t^3}$$

Thus, for some placement of the grids, the resulting $g_i$ satisfy
$$\sum_{i=1}^r (\langle H, g_i \rangle)^2 \ge \frac{\rho^2 \delta^{13} m^2
(\log n)^{1/2}}{2^{27} c^2_0 t^3}$$ yielding the statement of the main
lemma.  \qed \\

\noindent{\bf {Acknowledgement}} \\

I thank Professor J\'{o}zsef Beck for pointers to literature, useful
discussions and constant encouragement. \\

\noindent{\bf {References}} \\

\noindent 1. K. F. Roth, {\em {On irregularities of distribution}}. 
Mathematika 1, 1954. \\

\noindent 2. K. F. Roth, {\em {Remark concerning integer sequences}}. Acta
Arithmetica 9, 1964. \\

\noindent 3. P. Erd\H{o}s and J. Spencer, {\em {Probabilistic Methods in
Combinatorics}}. Akad\'emiai Kiad\'o, Budapest, 1974. \\

\noindent 4. P. Erd\H{o}s, {\em {On the combinatorial problems which 
I would most like to see solved}}. Combinatorica 1, 1981. \\

\noindent 5. J. Beck, {\em {Roth's estimate of the discrepancy of
integer sequences is nearly sharp}}. Combinatorica 1, 1981. \\

\noindent 6. J. Beck, {\em {On irregularities of $\pm 1$-sequences}}. 
\"{O}sterreich. Akad. Wiss. Math.-Natur. Kl. Sitzungsber II 195, 1986. \\

\noindent 7. R. Hochberg, {\em {Discrepancy and Bandwidth}}. Ph.D. Thesis, 
Rutgers University, 1994. \\

\noindent 8. J. Matou\v{s}ek and J. Spencer, {\em {Discrepancy in
arithmetic progressions}}. Journal of the American Mathematical Society 9,
1996. \\

\noindent 9. D. Reimer, {\em {Five Coloring Theorems}}. Ph.D. 
Thesis, Rutgers University, 1997. 

\end{document}